# THE ELEMENTS OF AN ANALYSIS OF THE FUNCTIONS OF A SET

## A. A. BOSOV

**Abstract.** The variant of calculation of functions of set and their application is offered. In particular: the new measure of system of sets generalizing classical concept of a measure is entered; the variation of set that has allowed to construct a derivative of function of set on a measure is entered, and to receive necessary conditions of optimization in terms of function of set.

**Introduction**

The notion of coexistent indexes [1] was formulated by O.L Cauchy. These indexes are determined by different circumstances: geometrical, physical, etc. In other word, if there is some physical body, then we can juxtapose its mass, volume, surface area and so on. Due to the Mathematics these indexes can be considered as some functions of set. In reality only the coexistent indexes can be determined, and we can call them as Cauchy`s indexes and the way of their mathematical description as the function of set.

The aim of the work is to create elements of an analysis of the functions of a set.

**1. Sets and systems of sets, measure.**

Let set $\Omega$ consist of finite number of elements $\omega_i$, $i = \overline{1,n}$, and the system $\mathfrak{A}(\Omega)$ is various subsets of set $\Omega$. Let's determine function of set from elements A of $\mathfrak{A}(\Omega)$ in this way:

$$m(A) = \begin{cases} 0, if\ A = \varnothing; \\ \sum_{w \in A} H(\omega), if\ \dot{A} \neq \varnothing. \end{cases} \qquad (1)$$

In this definition the function $H(\omega) \neq 0$ for any $\omega \in \Omega$.

Let set $\Omega$ represent $\Omega = \{X : 0 \leq X \leq 1\}$, let's consider interval like (a, в), [a, в), (a, в], [a, в] as a system of a subset $\Omega$, where $0 \leq a \leq b \leq 1$. This system is defined by symbol S($\Omega$), and let $\mathfrak{A}(S(\Omega))$ denote various subsets of the system S($\Omega$). If A is element of $\mathfrak{A}(S(\Omega))$, then in general his structure is A={K,Q}, where K – continual set, Q – set of isolated points.

Let's consider a function of a set m(A), which is denoted on A$\in \mathfrak{A}(S(\Omega))$ and determined by formula:

$$m(A) = \begin{cases} 0, if\ \ A = \varnothing; \\ \mu_l(K) + c \sum_{x \in Q} H(x),\ if\ \ A \neq \varnothing. \end{cases} \qquad (2)$$

where $\mu_l(\cdot)$ – Lebesque`s measure, c - positive constant.

Let F be Cantor`s (Cantorian) set.

System S($\Omega$) is replenished with F and defined by S$^*$($\Omega$).

If $\mathfrak{A}(S^*(\Omega))$ is a system (class) of various of subset S$^*$($\Omega$), then structure of any elements from $\mathfrak{A}(S^*(\Omega))$ represents A={K, f, Q}, where K – continual set or a suite of continual set, f – subset of the set F or a suite of fractals, Q is the set of isolated points or a suite of sets of isolated points. Let element A of $\mathfrak{A}(S^*(\Omega))$ correspond number, which is defined by formula:

$$m(A) = \begin{cases} 0, & \text{if } A = \varnothing; \\ \alpha_1 \mu_l(K) + \alpha_2 \mu_c(f) + \alpha_3 \sum_{x \in Q} H(x), & \text{if } A \neq \varnothing. \end{cases} \quad (3)$$

where $\mu_c(f)$ - Caratheodory`s measure of the set f[2], and the set f will be called as fractal [3], $\alpha_1, \alpha_2, \alpha_3$ - coefficients, that make corresponding components to be of common dimensions.

All three types of determination of $m(A)$ have next properties:

p1. $m(A) \geq 0$; if $m(A) = 0$, then $A = \varnothing$
p2. $m(A \cup B) = m(A) + m(B) - m(A \cap B)$. (4)

Hereinafter the function $m(A)$ will be called as a measure of a system of sets, which is determined on class $\mathfrak{A}(S^*(\Omega))$. Let $\rho(A,B) = m(A \Delta B)$, where $\Delta$ - operation of a symmetrical difference of sets A and B. In this case there exist:

1. $\rho(A,B) \geq 0$; if $\rho(\dot{A}, \hat{A}) = 0$, then $\dot{A} = \hat{A}$;
2. $\rho(\dot{A}, \hat{A}) = \rho(\hat{A}, \dot{A})$;
3. $\rho(A,B) \leq \rho(A,C) + \rho(C,B)$,

where $A, B, C \in \mathfrak{A}(S^*(\Omega))$.

So, $\rho(A,B)$ is the distance among elements of $\mathfrak{A}(S^*(\Omega))$, therefore, the pair $\langle \mathfrak{A}(S^*(\Omega)), \rho(\cdot, \cdot) \rangle$ is non-linear metric space, which is denoted by $\mathfrak{F}$. We need operation of variation of some set A in order to build an informal analysis of the functions of a set.

**Definition 1.** The set $C = A \Delta B$ is called as variation of set A by B.

Notice, that if $\mathfrak{A}(S^*(\Omega))$ is algebra, then $\forall A, B \in \mathfrak{A}(S^*(\Omega))$ and $A \Delta B \in \mathfrak{A}(S^*(\Omega))$.

## 2. The operation of limit.

Let sequence $\{B_n\}$, n=1,2,3,.... from $\mathfrak{A}(S^*(\Omega))$ be such a type of sequence that we can point out some set $B \in \mathfrak{A}(S^*(\Omega))$, which has next property: there is such a number $n(\varepsilon)$, $\varepsilon > 0$,, and $m(B_n \Delta B) < \varepsilon$ exists for $n > n(\varepsilon)$.

The set B is called as limit of sequence of sets $\{B_n\}$.

Notice, that if sequence $\{B_n\}$ has one more limit set $B'$, then $B = B'$. Really,

$$m(B' \Delta B) \leq m(B_n \Delta B') + m(B_n \Delta B),$$

so, from definition of limited set we can say that

$$m(B' \Delta B) < 2\varepsilon.$$

Due to arbitrariness $\varepsilon$ we obtain $m(B'\Delta B)=0$, and due to property 1 we have $B'\Delta B = \varnothing$, i.e. $B = B'$.

Let $\underline{B}$ be such a kind of set, that $\forall x \in \underline{B}$ and we can point out the number n(x), when dependency $x \in B_n$ exists for all n> n(x).

Let $\overline{B}$ denote such a kind of set, that $\forall x \in \overline{B}$ follows, that $x$ belongs to infinite number of sets from sequence {B$_n$}.

As we know [4], the set $\underline{B}$ is called as lower limit and $\overline{B}$ - as upper limit of sequence {B$_n$}. In case, when $\underline{B} = \overline{B}$, let B$_*$ denote this set and we will call it as limit of sequence {B$_n$} or limited Borel`s set.

**Theorem 1.** If the sequence {B$_n$} has limits B$_*$ and B, then these limits are equal; moreover the existence of each of them follows from existence of another.

**Proof**. Let $x \in B_*$, then we can point out such a number n(x), that $x \in B_n$ for n>n(x). Assume, that $x \notin B$, then $m(B_n \Delta B) \geq m(\{x\})$ we get, that $m(B_n \Delta B)$ doesn't tend to zero, when $n \to \infty$. Therefore, $x \in B$ and we have $B_* \subseteq B$. In this way we can prove, if $x \in B$, then $x \in B_*$ and then $B \subseteq B_*$, that proves the first part of the theorem.

Let $B_*$ exist. Prove, that $m(B_n \Delta B_*) \to 0$ if $n \to \infty$. Assume, that it is not correct, i.e. there is $y \in B_n$, but $y \notin B_*$, but $y$ belongs to infinite number of sets from the sequence $\{B_n\}$, then $y \in \overline{B}$, and there we find contradiction because $\overline{B} = B_*$. Therefore
$$\lim_{n\to\infty} m(B_n \Delta B_*) = 0 \to B_* = B.$$

We will use next statement in order to prove the existence of $B_*$, provided that there is $B$.

**Statement.** If sequence $\{B_n\}$ converges, then any of its subsequence also converges and they have the only one common limit.

Really, if $a_n = m(B_n \Delta B)$ then we'll get numeral sequence, that converges to zero and if $\{a_{n_k}\}$ is subsequence of sequence $\{a_n\}$, then $\lim_{n_k \to \infty} a_{n_k} = 0$ and because $a_{n_k} = m(B_{n_k} \Delta B)$ we get proof of statement.

Let $B$ exist, then $\forall x \in B$ it follows, that $x$ belongs to infinite number of sets from $\{B_n\}$ and that means that $x \in \overline{B}$. Let $x \notin \{B_{n_k}\}$, then $\lim_{n_k \to \infty} m(B_{n_k} \Delta B) \neq 0$, and that contradicts existence of $B$; therefore, $x$ belongs to all $B_n$ starting from certain number, i.e. $x \in \underline{B}$. If $\overline{B} \Delta \underline{B} \neq \varnothing$, then we have contradiction the existence $B$, and that is why $\overline{B} = \underline{B}$, that proves existence $B_*$.
Using accepted notation, we can write

$$\underline{\lim_{n\to\infty}} B_n = \underline{B}; \quad \overline{\lim_{n\to\infty}} B_n = \overline{B}; \quad \lim_{n\to\infty} B_n = B;$$

the two first limits always exist and mean

$$\underline{\lim_{n\to\infty}} B_n = \bigcup_{k=1}^{\infty} \bigcap_{n=k}^{\infty} B_n; \quad \overline{\lim_{n\to\infty}} B_n = \bigcap_{k=1}^{\infty} \bigcup_{n=k}^{\infty} B_n.$$

If $\underline{B} = \overline{B}$, then the Borel`s limit $B_* = \underline{B} = \overline{B} = B$ due to theorem 1. Therefore, if there is B, then there exists

$$\lim_{n\to\infty} B_n = \bigcup_{k=1}^{\infty}\bigcap_{n=k}^{\infty} B_n = \bigcap_{k=1}^{\infty}\bigcup_{n=k}^{\infty} B_n. \qquad (5)$$

Dependency (5) should be understood as representation of operation of limit by the operation of union and intersection of sets. On the other hand, $B_n$ are elements of the metric space $\langle \mathfrak{A}(S^*(\Omega)), \rho(\cdot,\cdot) \rangle$ and dependency (5) denotes operation of limit in this space.

**Theorem 2.** Let sequence $\{B_n\}$ have its limit B, then sequence $\{B_n * A\}$ has its limit $B * A$, where $*$ - the operation among $\cup, \cap, \setminus, \Delta$.

**Proof**. Let`s consider sequence $\{B_n \cup A\}$. Then $m((B_n \cup A)\Delta(B\Delta A)) = m(B_n \Delta B)$, and we will get

$$\lim_{n\to\infty} m((B_n \cup A)\Delta(B \cup A)) = \lim_{n\to\infty} m(B_n \Delta B) = 0;$$

Therefore, $B \cup A$ is limit of the sequence.

For the rest operations the proof is based on dependency $(B_n * A)\Delta(B * A) \subset B_n \Delta B$ and on monotony of measure.

### 3. Differentiation of the functions of a set [5]

Let F(A) be function of set, that is determined on algebra $\mathfrak{A}(\Omega)$, and m(A) is a measure, that is defined by one of the formula (1)-(3).

If $\{B_n\}$ is certain sequence of elements of $\mathfrak{A}(\Omega)$, then let's consider sequence of number, that can be determined by formula:

$$a_n = \frac{F(A\Delta B_n) - F(A)}{m(A\Delta B_n) - m(A)}; \quad n=1,2..... \qquad (6)$$

**Definition 2.** If sequence of numbers (6) has limit, then we'll call it as a derivative of a function of the set F(A) by the measure m(·) on the sequence $\{B_n\}$ and denote it in next way:

$$\left.\frac{dF(A)}{dm}\right|_{\{Bn\}} \stackrel{def}{=} \lim_{n\to\infty} a_n, \qquad (7)$$

Evidently, that numbers

$$\left.\underline{DF}(A)\right|_{\{B_n\}} = \varliminf_{n\to\infty} a_n; \quad \left.\overline{DF}(A)\right|_{\{B_n\}} = \varlimsup_{n\to\infty} a_n, \qquad (8)$$

always exist, and if they are finite and equal each other, then it is necessary and enough for existence of a derivative (7).

**Theorem 3.** If $F(E) \leq F(A), \quad \forall A \in \mathfrak{A}(S(\Omega))$, , then due to the necessary condition there is

$$sign\left(\left.\underline{DF}(E)\right|_{\{B_n\}\to B\subset E}\right) = sign\left(\left.\overline{DF}(E)\right|_{\{B_n\}\to B\subset E}\right) \leq 0,$$

and if there exists derivative (7), then

$$\left.\frac{dF(E)}{dm}\right|_{\{B_n\}\to B\subset E} \leq 0. \tag{9}$$

**Proof.** Notice, that $F(E\Delta B_n)-F(E)\geq 0$, and $m(E\Delta B_n)-m(E)=m(B_n)-2m(E\cap B_n)$..
Because $B\subset E$, we can point out such a number $n_*$, starting with which $B_n \subset E$.

According to the note, we will get, that if $n>n_*$, when A=E, then numbers $a_n$ in the sequence (6) will be determined by formula

$$a_n = \frac{F(E\Delta B_n)-F(E)}{-m(B_n)}$$

because if $m(B_n)>0$, then $a_n \leq 0$, and that proves theorem.

If limited set $B \not\subset E$, then, when there is derivative on sequence $\{B_n\}$, that converges to B, we'll get

$$\left.\frac{dF(E)}{dm}\right|_{\{B_n\}\to B\cap E=\varnothing} \geq 0 \quad . \tag{10}$$

In general case derivatives considered depend on E, sequence $\{B_n\}$ and limited set B.

We'll take class of continual functions of set for investigation of these dependencies.

**Definition 3.** Function of set F(A) will be called as continual, if we can point out such $\delta(A,\varepsilon)$, for any $\varepsilon>0$, that if only $m(A\Delta B)<\delta(A,\varepsilon)$, then it will cause

$$|F(B)-F(A)|<\varepsilon \ .$$

In case, when $\delta(A,\varepsilon)$ doesn't depend on A, function F(A) will be called as staedy continual on $\mathfrak{A}(S^*(\Omega))$.

**Theorem 4.** If F(A) is continual and sequence $\{B_n\}$ converges to B, then there exists

$$\lim_{n\to\infty} F(B_n) = F(\lim_{n\to\infty} B_n) \quad . \tag{11}$$

**Proof.** Due to convergence of sequence $\{B_n\}$ and continuity of F(A), when $\forall \varepsilon$, we can point out such $n(\varepsilon)$, that when $n>n(\varepsilon)$ we'll get

$$|F(B_n)-F(B)|<\varepsilon \ ,$$

i.e. numeral sequence $\{F(B_n)\}$ has its limit as number F(B), that proves formula (11).

**4. Rules of differentiation.**

There are next rules of differentiation of a function of a set by measure on the sequence $\{B_n\}$:

r.1 If $F(A)=const$, $\forall A \in \mathfrak{A}(S^*(\Omega))$, then derivative is equal 0;

r.2 Constant multiplier can be out of the symbol of derivative;

r.3 Derivative of sum of two functions is equal sum of derivatives, if the last exist;

in order to formulate next rules of differentiation, let's assume that the function is continual and it always has derivative on converging sequence, and the derivative equals

$$\left.\frac{dF(A)}{dm}\right|_{\{B_n\}\to B\neq\varnothing} = \frac{F(A\Delta B)-F(A)}{m(A\Delta B)-m(A)} \tag{12}$$

r.4 Derivative of product:

$$\frac{d(F_1(A)F_2(A))}{dm} = \frac{F_2(A\Delta B)+F_2(A)}{2}\cdot\frac{dF_1}{dm}+\frac{F_1(A\Delta B)+F_1(A)}{2}\cdot\frac{dF_2}{dm}$$

r.5 Derivative of fraction:

$$\frac{d(F_1(A)/F_2(A))}{dm} = \frac{F_2(A)\cdot\dfrac{dF_1}{dm}-F_1(A)\cdot\dfrac{dF_2}{dm}}{F_2(A\Delta B)\cdot F_2(A)}$$

r.6. $\dfrac{d\varphi(F(A))}{dm} = \dfrac{d\varphi(t)}{dt}\bigg|_{t=(F(A)+\theta(F(A\Delta B)-F(A)))}\dfrac{dF(A)}{dm}$

We need next lemma in order to prove the formula (12).

**Lemma 1.** If the sequence $\{B_n\}$ converges to B, and function $F_1(A)$ and $F_2(A)$ are continual, then $F_1(A)/F_2(A)$ is a continual function, when $F_2(A)\neq 0$ and

$$\lim_{n\to\infty}\frac{F_1(B_n)}{F_2(B_n)} = \frac{F_1\left(\lim_{n\to\infty}B_n\right)}{F_2\left(\lim_{n\to\infty}B_n\right)} \qquad (13)$$

**Proof.** Continuity of dependency is defined by inequality

$$\left|\frac{F_1(B)}{F_2(B)}-\frac{F_1(A)}{F_2(A)}\right| < \frac{|F_1(A)|+|F_2(A)|}{|F_2(B)\cdot F_2(A)|}\cdot\varepsilon,$$

where $\varepsilon$ is arbitrary and formula (13) is determined by theorem 4.

**Lemma 2.** $\lim_{n\to\infty} m(B_n) = m\left(\lim_{n\to\infty} B_n\right)$

**Proof.** This lemma states, that measure is a continual function of a set. Really, as $B_n\to B$, then $m(B\Delta B_n) = m(B\setminus B_n)+m(B_n\setminus B) < \varepsilon$ starting with number $n(\varepsilon)$. Let $\varepsilon_1 = m(B\setminus B_n)$; $\varepsilon_2 = m(B_n\setminus B)$, we have $\varepsilon_1+\varepsilon_2 < \varepsilon$, and considering that if

$$m(B) = \varepsilon_1 + m(B\cap B_n);$$
$$m(B_n) = \varepsilon_2 + m(B\cap B_n),$$

then

$$|m(B)-m(B_n)| = |\varepsilon_1-\varepsilon_2| < \varepsilon$$

and we'll get proof of the lemma due to arbitrariness of $\varepsilon$.

For proving the formula (12) we'll use lemma 1, lemma 2, and theorem 2. Formal note of proving is

$$\lim_{n\to\infty}\frac{F(A\Delta B_n)-F(A)}{m(A\Delta B_n)-m(A)} = \frac{\lim_{n\to\infty}(F(A\Delta B_n)-F(A))}{\lim_{n\to\infty}(m(A\Delta B_n)-m(A))} = \frac{F\left(A\Delta\lim_{n\to\infty}B_n\right)-F(A)}{m\left(A\lim_{n\to\infty}B_n\right)-m(A)}.$$

**5. Some examples of tasks of function of set.**

If set A is $\{\omega_1, \omega_2, ...\omega_k\}$, then elements' order in A isn't essential, and then function of set F(A) is considered as function of elements $\omega_1, \omega_2, ...\omega_n$, and it won't change its value whatever exchange may be.

Therefore, $F(\{\omega_1, \omega_2, ...\omega_n\})$ is symmetrical function of enumerated elements [6].

Let set $\Omega_n = \{\omega_1, \omega_2, ...\omega_n\}$ and number $x_i = x(\omega_i)$, $i = \overline{1, n}$, is juxtaposed for each element. If we determine elementary symmetrical functions

$$\sigma_1 = x_1 + x_2 + ... + x_n;$$
$$\sigma_2 = x_1 x_2 + x_1 x_3 + ... + x_{n-1} x_n;$$
$$\sigma_3 = x_1 x_2 x_3 + x_1 x_2 x_4 + ... + x_{n-2} x_{n-1} x_n;$$
$$\cdot \cdot \cdot \cdot \cdot \cdot \cdot \cdot \cdot \cdot \cdot$$
$$\sigma_n = x_1 x_2 ... x_n,$$

then any function $\varphi(\sigma_1, \sigma_2, ...\sigma_n)$ as function $x_1, x_2, ...x_n$ is symmetrical function.

If $\varphi(t_1, t_2, ...t_m)$ is arbitrary function of variables $t_i, i = \overline{1, m}$, then the function of set can be determined by next formula

$$F(A) = \varphi(\sigma_1(A), \sigma_2(A), ...\sigma_m(A)),$$

where m<n.

Notice, if m<n, then we can determine other functions of set by function $\varphi$, choosing any m of elementary functions $\sigma$ and locating them due to their place in $\varphi$ in different way. So, we can build n!/(n-m)! of various functions of set by m-ary function $\varphi$ for set A.

**Example 1**. Let $A = \{\omega_1, \omega_2, \omega_3\}$ and each element is juxtaposed with numbers $x_1, x_2, x_3$. If $\varphi = t_1^2 + t_2^3$, then we'll have n=3, m=2 and we can build 3!/(3-2)!=3!=6 of various functions of set, that are:

$F_1(A) = \sigma_1^2 + \sigma_2^3$; $F_2(A) = \sigma_2^2 + \sigma_1^3$; $F_3(A) = \sigma_1^2 + \sigma_3^3$; $F_4(A) = \sigma_3^2 + \sigma_1^3$;
$F_5(A) = \sigma_2^2 + \sigma_3^3$; $F_6(A) = \sigma_3^2 + \sigma_1^3$, ãäå $\sigma_1 = x_1 + x_2 + x_3$; $\sigma_2 = x_1 x_2 + x_1 x_3 + x_2 x_3$; $\sigma_3 = x_1 \cdot x_2 \cdot x_3$.

To sum it up, we can assume that next variant of building of function of set, if there is a suite of functions, $\varphi_i, i = \overline{1, n}$, where $\varphi_i$ is $i$-ary function.

Then

$$F(A) = \sum_{i=1}^{n} \varphi_i(\sigma_{i1}(A), \sigma_{i2}(A), ...\sigma_{ij}(A)),$$

that is function of set.

The number of such functions equals:

$$N = \sum_{i=1}^{n} \frac{n!}{(n-1)!}.$$

Let function of set F(A) be determined on $\mathfrak{A}(S(\Omega))$, where set $\Omega_n = \{\omega_1, \omega_2, ...\omega_n\} \subset \Omega$, then if $x_i = F(\{\omega_i\})$ and there are elementary symmetrical functions $\sigma_i(A)$, where $A \in \mathfrak{A}(S(\Omega))$, we'll consider function

$$\widetilde{F}_n(A) = \sum_{k=1}^{n} c_k \sigma_k(A), \tag{14}$$

which will be called as decomposition of base function F(A) by elementary symmetrical functions.

If we define scalar product by formula

$$\langle F_1(A), F_2(A) \rangle = \sum_{A \in (\Omega_n)} F_1(A) \cdot F_2(A),$$

then coefficient $c_k$ can be determined from system:

$$\langle F(A), \sigma_1(A) \rangle = c_1 \langle \sigma_1, \sigma_1 \rangle + c_2 \langle \sigma_2, \sigma_1 \rangle + \ldots + c_n \langle \sigma_n, \sigma_1 \rangle;$$
$$\langle F(A), \sigma_2(A) \rangle = \qquad\qquad c_2 \langle \sigma_2, \sigma_2 \rangle + \ldots + c_n \langle \sigma_n, \sigma_2 \rangle;$$
$$\cdot \qquad \cdot \qquad \cdot \qquad \cdot \qquad \cdot \qquad \cdot \qquad \cdot$$
$$\langle F(A), \sigma_n(A) \rangle = c_n \langle \sigma_n, \sigma_n \rangle;$$

In this case, function $\widetilde{F}(A)$ in regard to function F(A) is such a function, that is minimum

$$S^2(c_1, c_2, \ldots c_n) = \sum_{A \in (\Omega_n)} \left( F(A) - \sum_{k=1}^{n} c_k \sigma_k(A) \right)^2,$$

i.e. $\widetilde{F}(A)$ is approximation of function F(A) due to the method of least squares on $\mathfrak{A}(S(\Omega))$. Notice, that if $c_0 = F(\varnothing)$ and if you define unit function for any $A \in \mathfrak{A}(S(\Omega))$, then representation

$$\widetilde{F}_n(A) = \sum_{k=0}^{n} c_k \sigma_k(A)$$

can be considered as Fourier n-series

$$\widetilde{F}(A) = \sum_{k=0}^{\infty} c_k \sigma_k(A)$$

for function F(A) on the set $(\Omega_\infty)$, where $\Omega_\infty$ is countable set, that belongs to base set $\Omega$.

## 6. Integral by measure $m(\bullet)$.

Let
$$f(x) = F(\{x\})$$

and if $\alpha \le f(x) \le \beta$, when $x \in A$, then we'll split interval $[\alpha, \beta]$ in this way:

$$\alpha = y_0 < y_1 < \ldots < y_n = \beta$$

and let
$$e_\nu = \{x \in A : y_\nu \le f(x) < y_{\nu+1}\}.$$

Using Lebeque`s scheme, we'll get

$$\underline{S}_n = \sum_{\nu=0}^{n} y_\nu m(e_\nu),$$

$$\overline{S}_n = \sum_{\nu=0}^{n} y_{\nu+1} m(e_\nu),$$

where $y_{n+1} = \beta$.

If
$$\lim_{n \to \infty} \underline{S}_n = \lim_{n \to \infty} \overline{S}_n = J(A)$$

when $\max(y_{\nu+1} - y_\nu) \to 0$, then common limit will be called as Lebeque`s integral by measure $m(\bullet)$, and will be denoted in this way:

$$J(A) = \int_A f(x) dm(x).$$

The problem of existence of this integral is similar to the situation, when it is integrated by Lebeque`s measure [1].

Notice, if set A consists of discrete number of points $A=\{x_1, x_2, \ldots x_n, \ldots\}$, then the integral can be determined by formula

$$J(A) = c \sum_{x_i \in A} f(x_i) H(x_i).$$

In case when A consists of interval $[a,b] \in [0,1]$ and set $Q=\{x_1, x_2, \ldots x_n, \ldots\}$, then the integral can be determined by formula

$$J(A) = \int_a^b f(x) d\mu(x) + c \sum_{x_i \in Q} f(x_i) H(x_i),$$

where $\mu$ is usual Lebeque`s measure.

In the formula some $x_i \in Q$ can belong to interval [a,b].

So, if A doesn't include set Q, then integral defined will coincide with classical Lebeque`s integral.

The function H(x) is a measure of point x and in every certain case its representation is determined by content integral value, and randomness in its representation let it be chosen, so that mathing series will converge.

We can show the last statement on an example, when f(x) is limited on real axis, and set $Q=\{\ldots -n, -(n-1) \ldots 0, 1, 2, \ldots n, ..\}$, then take $H(x) = \exp(-|x|)$. Convergence of series $\sum_{x \in Q} f(x) H(x)$ is evident and by that the existence of integral defined is evident too.

If A is such a set, when $m(A) < \infty$, and function f(x) is limited on A, then there is next theorem.

**Theorem 5**. (mean- value thearem)

$$\alpha m(A) \leq \int f(x) dm(x) \leq \beta m(A),$$

where $\alpha = \inf f(x), \beta = \sup f(x)..$

On the base of the theorem we we'll get Radon-Nikodim`s theorem

$$\left. \frac{dJ(A)}{dm} \right|_{\{B_n\} \to \{x\}} = f(x),$$

if f(x)-continual function.

**7. n-ary functions of set.**

**Definition 4.** Let $\Theta = \overset{n}{\underset{i=1}{X}} \mathfrak{A}(\Omega_i)$ be Cartesian product of classes $\mathfrak{A}(\Omega_i)$, then reflection $\Theta$ on real axis will be called as n-ary function of set and defined by note

$$F(A_1, A_2, \ldots A_n),$$

where the list $[A_1, A_2, \ldots A_n] \in \Theta$.

**Definition 5.** The limit

$$\lim_{n \to \infty} \frac{F(A_1, A_2, \ldots A_i \Delta B_n^i, \ldots A_n) - F(A_1, A_2, \ldots A_i, \ldots A_n)}{m_i(A_i \Delta B_n^i) - m(A_i)},$$

if it exists, will be called as a partial derivative by measure $\mu_i$ on sequence $\{B_n^i\}$ and denoted by:

$$\left.\frac{\partial F}{\partial m_i}\right|_{\{B_n^i\}}.$$

**Definition 6.** Vector $\nabla F$, whose components are equal $\left.\dfrac{\partial F}{\partial m_i}\right|_{\{B_n^i\}}$ i=1..n , will be called as gradient of function $F(A_1, A_2, \ldots A_n)$.

**Theorem 6.** If $F(A_1, A_2, \ldots A_n)$ -min , then due to the necessary condition there is

$$\nabla F^-(A_1, A_2, A_3) \leq 0,$$

where $\nabla F^-(A_1, A_2, A_3)$ is inner gradient, whose components are

$$\left.\frac{\partial F}{\partial m_i}\right|_{\{B_n^i\} \to B^i \subseteq A_i}.$$

The theorem is generalization of theorem 3 on n-ary function of set.

**8. Some tasks of vector optimization.**
**Example 2**.
  There are two factors, that are considered:

$$F_1(A) = \int_A f_1(x) dm;$$

$$F_2 = \int_{\Omega \setminus A} f_2(x) dm,$$

and each one should be minimized.
  Officially, the task should be defined as

$$\begin{pmatrix} F_1(A) \\ F_2(A) \end{pmatrix} \to \min,$$

provided that $A \subset \mathfrak{A}(S^*(\Omega))$ .
  Notice, the solution of the task is considered to be next:

$$\mathfrak{A}^*(S^*(\Omega)) \subseteq \mathfrak{A}(S^*(\Omega))$$

where $\forall A, B \subset \mathfrak{A}^*(S^*(\Omega))$, there are

$$F_1(A) \leq F_1(B);$$
$$F_2(A) \geq F_2(B),$$

or opposite inequalities. The variant of equalities isn't excepted.
If $A$ is one of solution, then we`ll obtain

$$\frac{dF_1(A)}{dm} = -\lambda \frac{dF_2(A)}{dm},$$

where $\lambda \geq 0$.
For considering functions $F_1, F_2$ and sequence $\{B_n\}$, that converges to $x \in A$, the necessary condition is:

$$f_1(x) = \lambda f_2(x).$$

If $\Omega = \{(x_1, x_2) : x_1, x_2 \geq 0, x_1^2 + x_2^2 \leq 1\}$, $f_1 = x_1 + x_2, f_2 = ax_1 - bx_2, a, b > 0$, then we'll get

$$x_1 + x_2 = \lambda(ax_1 - bx_2),$$

or

$$x_2 = \frac{\lambda a - 1}{\lambda b + 1} x_1,$$

where $\lambda a - 1 \geq 0$.
So, when $\lambda \geq 0$, $\lambda$ is fixed, the set $A(\lambda) \subset \mathfrak{A}^*(S^*(\Omega))$ is:

$$A(\lambda) = \bigcup_{\mu=0}^{\lambda} \{(x_1, x_2) : x_2 = \frac{\lambda a - 1}{\lambda b + 1} x_1\} \bigcap \Omega.$$

**Example 3.** Let n-ary function of sets be

$$F(\Omega_1, \Omega_2, \Omega_3) = \sum_{i=1}^{3} \int_{\Omega_i} f_i(x) dm,$$

which is needed to be minimized, provided that

$$\bigcup_{i=1}^{3} \Omega_i = \Omega; \quad \Omega_i \bigcap \Omega_j = \emptyset, i \neq j.$$

Due to the admissible variation and theorem 6, the task solution is:

$$\Omega_1 = \{x \in \Omega : f_1(x) \leq f_2(x) \text{ and } f_1(x) \leq f_3(x)\};$$
$$\Omega_2 = \{x \in \Omega : f_2(x) < f_1(x) \text{ and } f_2(x) \leq f_3(x)\};$$
$$\Omega_3 = \{x \in \Omega : f_3(x) < f_1(x) \text{ and } f_3(x) < f_2(x)\}.$$

**8. Lagrange`s formula for functions of set.**
Let the function of set $F(A)$ be, $\forall A \in \mathfrak{A}(\Omega)$, and $m(A)$ is measure on $\mathfrak{A}(\Omega)$.

Let's build curve $\gamma$ on coordinates (x,y) in this way:

$$\gamma := \{(x,y) : x = m(A), y = F(A), A \in \mathfrak{A}(\Omega)\} .$$

If $y = \varphi(x)$ is dependency, that describes curve $\gamma$, we'll get derivative

$$\frac{F(B) - F(A)}{m(B) - m(A)} = \left. \frac{d\varphi(x)}{dx} \right|_{x = m(A) + \theta(m(B) - m(A))} .$$

where $\theta \in [0,1]$ .

Evidently, that the dependency is Lagrange's formula for functions of set.
Notice, that if function $F(A)$ is continual, then the function $\varphi(x)$ is also continual.
Dependency between derivatives of these functions is denoted by next formula

$$\lim_{B \to \varnothing} \frac{F(A \Delta B) - F(A)}{m(A \Delta B) - m(A)} = \left. \frac{d\varphi(x)}{dx} \right|_{x = m(A)} .$$